\documentclass[a4paper,10pt]{article}
\usepackage{amsfonts}
\usepackage{amstext}
\usepackage{amsthm}
\usepackage{amsmath}
\usepackage{amssymb}
\usepackage{latexsym}
\usepackage[latin1]{inputenc}
\newcommand{\bl}{\hfill\rule{2mm}{2mm}}
\newcommand{\R}{\mathbb{R}}

\newtheorem{mainthm}{Theorem}

\newcommand{\n}{\noindent}

\begin{document}

\title{Sharp $L^p$-Moser inequality on Riemannian manifolds
 }
\author{\textbf{Marcos Teixeira Alves \footnote{\textit{e-mail addresses}: mtmarcos@gmail.com(M. T. Alves)}}; \; \textbf{Jurandir Ceccon \footnote{\textit{E-mail addresses}:
ceccon@ufpr.br (J. Ceccon)}}\\
{\small\it Departamento de
Matem\'{a}tica, Universidade Federal do Paran\'{a},}\\
{\small\it Caixa Postal 019081, 81531-990, Curitiba, PR,
Brazil}}\maketitle

\markboth{abstract}{abstract}
\addcontentsline{toc}{chapter}{abstract}

\hrule \vspace{0,2cm}

\n {\bf Abstract}

We consider $(M,g)$ a smooth compact Riemannian manifold of dimension
$n \geq 2$ without boundary, $1 < p$ a real parameter and $r = \frac{p(n + p)}{n}$. This paper concerns the validity of the optimal Moser inequality

\[
\left( \int_M |u|^r\; dv_g \right)^{\frac{\tau}{p}} \leq \left(
A(p,n)^{\frac{\tau}{p}} \left(\int_M |\nabla_g u|^p\; dv_g\right)^{\frac{\tau}{p}}
+ B_{opt} \left(\int_M |u|^p\; dv_g\right)^{\frac{\tau}{p}}
\right) \left( \int_M |u|^p\; dv_g \right)^{\frac{\tau}{n}} \; .
\]

This kind of inequality was already studied in the last years in the
particular cases $1 < p < n$. Here we solve the case $n \leq
p$ and we introduce one more parameter $1 \leq \tau \leq
\min\{p,2\}$.  Moreover, we prove the existence of an extremal
function for the optimal inequality above.

\vspace{0,5cm} \hrule\vspace{0.2cm}

\section{Introduction}

Optimal inequalities of Moser \cite{Mo} type have been extensively
studied both in the Euclidean and Riemannian contexts. We refer the reader to
\cite{Ma}, \cite{Ba}, \cite{Beck}, \cite{CoVi}, \cite{DPDo} for
the Euclidean case and to \cite{Bro}, \cite{CMMZ}, \cite{CMJD},
\cite{CS} for Riemannian manifolds. In 1961, Moser \cite{Mo}
proved that solutions of certain elliptic equations of second
order have the standard $L^\infty$ norm dominated by the standard
$L^p$ norm for all $p > 1$. This technique has been improved and
is now known as Giorgi-Nash-Moser, which plays an important
role in the theory of PDEs. The idea developed by Moser to prove
that increase is based on an iteration process. The key point of
this technique in \cite{Mo} consists in connecting the solution a
particular PDE to inequality:

\[
\int_{\R^n} |u|^{\frac{2(n + 2)}{n}} dx \leq c \int_{\R^n} |\nabla u|^2 dx \left(\int_{\R^n} |u|^2 dx\right)^{\frac{2}{n}}
\]

\n where this inequality is valid for every function $u \in
C_0^\infty(\R^n)$ and some constant $c > 0$. Our interest herein
will be studying this inequality from the ``optimal viewpoint".

\n In recent years, this type of inequality has been explored
considering the best constant, \textit{i.e.}, the least possible constant 
in this inequality (we will make this precise below). In this
context, many authors are getting to apply this successfully in the
investigation of various problems in PDEs. On the famous paper of 
Brezis and Nirenberg \cite{brezis}, the key point to prove
existence of solution to nonlinear problems with critical exponent
was the optimal Sobolev inequality. Recently Del Pino and
Dolbeault \cite{DelDo} made use of the Gagliardo-Nirenberg optimal
inequality to study the optimal decay rate of the intermediate
asymptotics of solutions to nonlinear diffusion equations. The
optimal Gagliardo-Nirenberg inequality was also used by Farah
\cite{farah} to establish sufficient conditions for global
existence for Schrödinger equation. In particular, the optimal
inequalities of Moser type in Riemannian manifolds studied in
Ceccon and Montenegro \cite{CMMZ} (case $n = p = 2$) was applied
by Kishimoto and Maeda \cite{KM} to obtain global existence
theorems for Zakharov system in $\mathbb{T}^2$.

Denote by $W^{1,p}(\R^n)$ the classical Sobolev space. The \textit{general Euclidean Moser inequality} states that there exists $A > 0$ such that, for any function $u \in W^{1,p}(\R^n)$,

\begin{gather}\label{dgne}
\int_{\R^n} |u|^r\; dx \leq A
\left( \int_{\R^n} |\nabla u|^p\; dx \right) \left( \int_{\R^n}
|u|^p\; dx \right)^{\frac{p}{n}}\, ,
\tag{$M_E(A)$}
\end{gather}

\n where $1 \leq p$, $2 \leq n$ and $r = \frac{p(n + p)}{n}$. Define

\[
A(p,n)^{-1} = \inf_{u \in W^{1,p}(\R^n)} \{ ||\nabla
u||_{L^p(\R^n)}^p ||u||_{L^p(\R^n)}^{\frac{p^2}{n}}; ||u||_{L^r(\R^n)} = 1\} \; .
\]

\n The inequality $(M_E(A(p,n)))$ is called {\em optimal Euclidean Moser
inequality} and the constant $A(p,n)$ is the {\em best constant}
in this inequality. This optimal Moser constant was already studied by
Beckner in \cite{Beck}.

We now consider the Riemannian case. Let $(M,g)$  be a smooth
compact Riemannian manifold without boundary of dimension $2 \leq
n$. Using standard arguments (see \cite{DHV} for example) we
obtain a Riemannian version of the Euclidean inequality $(M_E(A))$.
Once there exists positive constants $C,D$ such that, for all
$u$ in the classical Riemannian Sobolev space $H^{1,p}(M)$, we
have

\begin{equation} \label{Riemann-non-sharp-no-tau}
\int_M |u|^r\; dv_g \leq \left( C \int_M |\nabla_g u|^p\; dv_g + D
\int_M |u|^p\; dv_g \right) \left( \int_M |u|^p\; dv_g
\right)^{\frac{p}{n}}\, ,
\end{equation}

\n where $1 \leq p$, $2 \leq n$ and $r = \frac{p(n + p)}{n}$.

We shall study a generalization of this inequality. Consider one more
parameter $\tau \in \R$, where $1 \leq \tau \leq p$. We immediately
have the inequality

\begin{gather}\label{AB1}
\left( \int_M |u|^r\; dv_g \right)^{\frac{\tau}{p}} \leq \left( A
\left(\int_M |\nabla_g u|^p\; dv_g\right)^{\frac{\tau}{p}} + B
\left(\int_M |u|^p\; dv_g\right)^{\frac{\tau}{p}} \right) \left(
\int_M |u|^p\; dv_g \right)^{\frac{\tau}{n}}\, , \tag{$M_R(A,B)$}
\end{gather}

\n for all $u \in H^{1,p}(M)$ and $p,n,r$ as above. Note that when $\tau = p$
we recover (\ref{Riemann-non-sharp-no-tau}). Observe that the non-sharp inequality  $(M_R(A,B))$ implies that

\begin{equation}\label{desin}
A \geq A(p,n)^{\frac{\tau}{p}} \; ,
\end{equation}

\n for any $1 \leq p$ and $2 \leq n$. This is shown by taking a
function with support contained in a small enough normal
neighborhood according to what was already observed in \cite{DHV}.

We now study the optimal inequality. Having two constants, the
optimality can be defined in two ways. We follow the more
interesting one from the PDE viewpoint (see chapters 4 and 5 of
the book \cite{He}). Define the {\em first Riemannian $L^p$-Moser
optimal constant} by

\[
A_{opt} = \inf \{ A \in \R:\; \mbox{there exists} \hspace{0,18cm} B
\in \R \hspace{0,18cm} \mbox{such that} \hspace{0,18 cm} M_R(A,B)
\hspace{0,18cm} \mbox{is valid}\}\, .
\]

\n Using a local argument and the definition of the optimal constant it is
easy shown

\[
A_{opt} \geq A(p,n)^{\frac{\tau}{p}} \; ,
\]

\n for any $1 \leq p$ and $2 \leq n$. Using a partition
of unit argument (the ideas of the outline are contained in \cite{DHV}) we
can establish that the first optimal Riemanniann $L^p$-Moser
inequality means that there exists a constant $C_\varepsilon \in
\R$ such that, for any $u \in H^{1,p}(M)$,

\[
\left( \int_M |u|^r\; dv_g \right)^{\frac{\tau}{p}} \leq \left(
(A(p,n)^{\frac{\tau}{p}} + \varepsilon) \left(\int_M |\nabla_g u|^p\;
dv_g\right)^{\frac{\tau}{p}} + C_\varepsilon \left(\int_M |u|^p\;
dv_g\right)^{\frac{\tau}{p}} \right) \left( \int_M |u|^p\; dv_g
\right)^{\frac{\tau}{n}}\, ,
\]

\n is valid for all $\varepsilon > 0$. Then, it is immediate that
the Euclidean optimal constant coincides with the first Riemannian
optimal constant, ie

\[
A_{opt} = A(p,n)^{\frac{\tau}{p}} \; .
\]

\n Since $(M_R(A(p,n)^{\frac{\tau}{p}} + \varepsilon,C_\varepsilon))$ holds, we can
define for all $\varepsilon > 0$:

\[
B_\varepsilon = \inf\{B \in \R; M_R(A(p,n)^{\frac{\tau}{p}} + \varepsilon,B)
\hspace{0,2cm} \mbox{is} \hspace{0,2cm} \mbox{valid}\} \; .
\]

\n Since constant non-zero functions belong to $H^{1,p}(M)$  the
constant $B_\varepsilon$ satisfies

\begin{equation}\label{sco}
B_\varepsilon \geq |M|^{- \frac{\tau}{n}} \; ,
\end{equation}

\n where  $|M|$ denotes the volume of $(M,g)$.

Then, we mean that for all $u \in H^{1,p}(M)$,

\begin{equation}\label{epsilon}
\left( \int_M |u|^r\; dv_g \right)^{\frac{\tau}{p}} \leq \left(
(A(p,n)^{\frac{\tau}{p}} + \varepsilon) \left(\int_M |\nabla_g u|^p\;
dv_g\right)^{\frac{\tau}{p}} + B_\varepsilon \left(\int_M |u|^p\;
dv_g\right)^{\frac{\tau}{p}} \right) \left( \int_M |u|^p\; dv_g
\right)^{\frac{\tau}{n}}\, ,
\end{equation}

\n is valid. It is easy to see that for $\varepsilon_1 <
\varepsilon_2$, we have that $B_{\varepsilon_2} \leq
B_{\varepsilon_1}$, ie, $B_\varepsilon$ is monotonous
non-increasing. Define

\[
{\cal B} = \lim_{\varepsilon \rightarrow 0} B_\varepsilon \;.
\]

In contrast with the Euclidean case ($B_\varepsilon = 0$), the
validity of the optimal inequality is delicate, since as
$\varepsilon \to 0$ the corresponding $B_\varepsilon$ might, at first, 
go to infinity. In fact, when $\tau = p > 2$  there
exist cases where the optimal inequality for Gagliardo-Nirenberg
or Sobolev is not valid,  depending on the geometry of $(M,g)$,
(see \cite {CMMZ} or \cite{D5} respectively). In this paper we
consider the parameter $\tau \leq \min\{p,2\}$, \textit{i.e.}, we weaken the
Moser inequality. This idea had already been applied in \cite{CS}
for a family of Gagliardo-Nirenberg inequalities when $\tau = 2$ and
$2 < p <n$. The ideas used in our proof were already present
in \cite{CS}, whose key point was to use techniques of explosion.
However, we introduce a new way to study this explosion
considering wich is more appropriate for this rescheduling.

\n The main objective of this paper is to show that ${\cal B}$ is
finite.

We now state the main results of this paper:

\begin{mainthm}\label{tgno1}
Let $(M,g)$ be a smooth compact Riemannian manifold without
boundary of dimension $n \geq 2$, $1 < p$
and $r = \frac{p(n + p)}{n}$. If $1 \leq \tau \leq \min\{p,2\}$
then ${\cal B} < \infty$ and $M_R(A(p,n)^{\frac{\tau}{p}},{\cal B})$ is always valid for all $u \in H^{1,p}(M)$.
\end{mainthm}

\n The Theorem \ref{tgno1} extends the results from Brouttelande
\cite{Bro} when $p = q =\tau = 2$, Ceccon and Montenegro
\cite{CMJD} in the case $1 < p = q \leq 2$, $\tau = p$ and Chen
and Sun \cite{CS} when $2 < p = q < n$ and $\tau = 2$.

The Theorem \ref{tgno1} allows to consider the {\em the second optimal constant}:

\[
B_{opt} = \inf \{B; M_R(A(p,n)^{\frac{\tau}{p}},B) \; \mbox{is valid} \} .
\]

\n It is immediate that $(M_R(A(p,n)^{\frac{\tau}{p}}, B_{opt}))$
is valid for all $u \in H^{1,p}(M)$. We call \textit{extremal function}
any non-zero function in $H^{1,p}(M)$ satisfying the equality
$(M_R(A(p,n)^{\frac{\tau}{p}}, B_{opt}))$. Unlike the Euclidean
case establish the existence of extremal functions is not
immediate just applying variational
techniques. But we could unify the study of the validity of the
optimal inequality with the existence of extremal functions. The
central point to establish the validity of the optimal inequality
and the existence of extremal function was to show that ${\cal B}
< \infty$. Thus, as a consequence of Theorem \ref{tgno1}, we have

\begin{mainthm}\label{extremal}
Let $(M,g)$ be a smooth compact Riemannian manifold without
boundary of dimension $2 \leq n$, $1 < p$ and $r = \frac{p(n + p)}{n}$. If $1 \leq \tau < \min\{2,p\}$ then
$M_R(A(p,n)^{\frac{\tau}{p}},B_{opt})$ admits an extremal function and $B_{opt} =
{\cal B}$.
\end{mainthm}

The Moser inequality is closely related to the very important
inequality of entropy. They are related through the Jensen's inequality.
The optimal Riemannian entropy inequality has recently been
studied in the case $1 < p \leq 2$ in \cite{ceccon-montenegro-logarithmic}.
The case $p > 2$ is open.

We will obtain now an estimate for the first best constant of
entropy inequality. For this, we first establish the entropy
inequality from Riemannian Moser inequality. Using the Jensen inequality in
$(M_R(A(p,n)^{\frac{\tau}{p}},{\cal B}))$ (as the ideas contained in
\cite{Beck}) we find the Riemannian entropy inequality

\begin{gather}
\int_M |u|^p\ln|u|^p dv_g \leq \frac{n}{\tau} \ln \left(A(p,n)^{\frac{\tau}{p}} \left(\int_M |\nabla_g u|^p dv_g \right)^{\frac{\tau}{p}} + {\cal B}\right)
\tag{$Ent(A(p,n)^{\frac{\tau}{p}},{\cal B})$}
\end{gather}

\n for all $u \in H^{1,p}(M)$ such that $||u||_{L^p(M)} = 1$,
where $2 \leq n$ and $1 < p$.

By proceeding analogously the Riemannian Moser inequality, we
define the \textit{first best constant for entropy inequality}. Denote by
$E(p,\tau,n)$ this optimal constant. That is, if $(Ent(A,B))$ is
valid for all function $u \in H^{1,p}(M)$ with $||u||_{L^p(M)} =
1$ then $A \geq E(p,\tau,n)$.

\n Thus, there is a constant $C_\varepsilon$ such that

\[
\int_M |u|^p \ln|u|^p dv_g \leq \frac{n}{\tau} \ln
\left((E(p,\tau,n) + \varepsilon) \left(\int_M |\nabla_g u|^p dv_g
\right)^{\frac{\tau}{p}} + C_\varepsilon\right)
\]

\n for all $u \in H^{1,p}(M)$ such that $||u||_{L^p(M)} = 1$ and
$\varepsilon > 0$.

For all $\varepsilon > 0$, we will also define

\[
E_\varepsilon = \inf\{B; Ent(E(p,\tau,n) + \varepsilon,B) \;
\mbox{is valid}\} \; .
\]

\n We will have, for each $\varepsilon > 0$, that $(Ent(E(p,\tau,n)
+ \varepsilon, E_\varepsilon))$ is valid for all $u \in H^{1,p}(M)$
such that $||u||_{L^p(M)} = 1$, $n \geq 2$ and $p > 1$. We must
note that it is not evident that $\limsup E_\varepsilon <
\infty$.

By the definition of the best constant $E(p,\tau,n)$ and as
$(Ent(A(p,n)^{\frac{\tau}{p}}, {\cal B}))$ is valid we can notice
that $E(p,\tau,n) \leq A(p,n)^{\frac{\tau}{p}}$. Applying again
the Jensen's inequality in entropy inequality $(Ent(E(p,\tau,n) +
\varepsilon,E_\varepsilon))$ (adapting the same ideas used in
\cite{ceccon}), we obtain the Nash inequality. Thus, we have
$E(p,\tau,n) \geq N(p,q,n)^{\frac{\tau}{p}}$, where $1 < q < p$
and $N(p,q,n)$ is the best constant in Nash inequality. This, we have the estimate

\[
N(p,q,n)^{\frac{\tau}{p}} \leq E(p,\tau,n) \leq A(p,n)^{\frac{\tau}{p}}
\]

\n  for all $q \in (1,p)$. As observed in \cite{ceccon}, we have
the limit

\[
\lim_{q \rightarrow p^-} N(p,q,n) = E_e(p,n)
\]

\n where $E_e(p,n)$ is the best constant in the Euclidean entropy
inequality obtained by Gentil \cite{gentil}. Then, we have the
following estimate for the first best contant for Riemannian
entropy

\[
E_e(p,n)^{\frac{\tau}{p}} \leq E(p,\tau,n) \leq A(p,n)^{\frac{\tau}{p}}
\]

\n for $2 \leq n$, $1 < p$ and $1 \leq \tau \leq \min\{p,2\}$.

For interesting applications of the entropy inequality see Gentil
\cite{gentil2} or Grillo \cite{grillo} and other references
contained therein.

\section{Proof of Theorem \ref{tgno1}}

To facilitate the understanding of the proof, we will divide the
argument in four steps.

- In section 2.1 we associate a sequence that satisfies an auxiliary minimal energy PDE;

- In section 2.2 we will prove that a phenomenon of explosion
occurs for this sequence;

- In section 2.3 we estimate the speed with which this sequence
converges to zero;

- In section 2.4 we will show that the constant $B_\varepsilon$ is
limited.

\subsection{The Euler-Lagrange equation associated with inequality of Moser} \label{EL}

\n As noted in the introduction (see inequality (\ref{epsilon}))
it suffices to show that

\[
\lim_{\varepsilon \rightarrow 0} B_\varepsilon = {\cal B} < \infty
\; .
\]

\n By (\ref{sco}) two possibilities might occur:

\vspace{0,2cm}

\n (C.1) ${\cal B} = |M|^{-\frac{\tau}{n}}$ or

\vspace{0,2cm}

\n (C.2) ${\cal B} > |M|^{-\frac{\tau}{n}}$.

\vspace{0,2cm}

\n If (C.1) occur, the Theorem 1 is proved by making $\varepsilon
\rightarrow 0$ in $M_R(A(p,n)^{\frac{\tau}{p}} + \varepsilon,
B_\varepsilon)$.

\n We will consider the case (C.2). Thus, there exists a
sequence $\gamma_\varepsilon > 0$ such that

\[
B_\varepsilon > |M|^{-\frac{\tau}{n}} + \gamma_\varepsilon \;.
\]

\n Just for simplicity, we consider that $\gamma_\varepsilon
\rightarrow 0$ when $\varepsilon \rightarrow 0$.

\n Now we will associate an Euler-Lagrange equation to the inequality
$(M_R(A(p,n)^{\frac{\tau}{p}} + \varepsilon, B_\varepsilon))$.

\n Define ${\cal H} = \{u \in H^{1,p}(M); \int_M |u|^r dv_g = 1\}$
and

\[
J_\varepsilon(u) = \left(A(p,n)^{\frac{\tau}{p}} \left(\int_M
|\nabla_g u|^p dv_g\right)^{\frac{\tau}{p}} + (B_\varepsilon -
\gamma_\varepsilon) \left(\int_M |u|^p dv_g
\right)^{\frac{\tau}{p}}\right)\left(\int_M |u|^p dv_g
\right)^{\frac{\tau}{n}} \; .
\]

\n Then, by definition of $B_\varepsilon$, there is a function
$u_0 \in {\cal H}$ such that

\[
J_\varepsilon(u_0) < 1 \; .
\]

\n Consider

\[
c_\varepsilon = \inf_{u \in {\cal H}} J_\varepsilon(u) < 1 \; .
\]

\n By a standard argument, it is easy to verify that $c_\varepsilon >
0$. In particular

\begin{equation}\label{um}
c_\varepsilon \left(\int_M |u|^r\; dv_g\right)^{\frac{\tau}{p}} \leq
\left( A(p,n)^{\frac{\tau}{p}} \left(\int_M |\nabla_g u|^p\;
dv_g\right)^{\frac{\tau}{p}} + (B_\varepsilon - \gamma_\varepsilon)
\left(\int_M |u|^p\; dv_g\right)^{\frac{\tau}{p}} \right) \left(
\int_M |u|^p\; dv_g \right)^{\frac{\tau}{n}}
\end{equation}

\n for all function $u \in H^{1,p}(M)$. Consider $u_k \in {\cal
H}$ such that $J_\varepsilon(u_k) \rightarrow c_\varepsilon$.
Let $\sigma > 1$ such that $c_\varepsilon < \sigma
c_\varepsilon < 1$. Then $J_\varepsilon(u_k) < \sigma
c_\varepsilon$ for $k$ sufficiently large. Using it together with
$M_R(A(p,n)^{\frac{\tau}{p}} + \lambda,B_\lambda)$ for $\lambda
> 0$ sufficiently small such that $A(p,n)^{\frac{\tau}{p}} -
\sigma c_\varepsilon(A(p,n)^{\frac{\tau}{p}} + \lambda) > 0$ and
as $p < r$, we can easily verify that the sequence $u_k$ is
bounded in $H^{1,p}(M)$. Therefore there is $\tilde{u}_\varepsilon
\in H^{1,p}(M)$ such that (for some subsequence) $u_k
\rightharpoonup \tilde{u}_\varepsilon$. Consequently

\[
u_k \rightarrow \tilde{u}_\varepsilon
\]

\n in $L^p(M) \cap L^r(M)$. In particular $||u_k||_{L^p}
\rightarrow ||\tilde{u}_\varepsilon||_{L^p}$ and $1 =
||u_k||_{L^r(M)} \rightarrow ||\tilde{u}_\varepsilon||_{L^r(M)}$.
So $\tilde{u}_\varepsilon \in {\cal H}$.

\n Therefore, we obtain

\[
J_\varepsilon(\tilde{u}_\varepsilon) \leq \lim_{k \rightarrow \infty}
\left(A(p,n)^{\frac{\tau}{p}} \left(\int_M |\nabla_g u_k|^p
dv_g\right)^{\frac{\tau}{p}} + (B_\varepsilon -
\gamma_\varepsilon) \left(\int_M |u_k|^p dv_g
\right)^{\frac{\tau}{p}} \right) \left(\int_M |u_k|^p dv_g
\right)^{\frac{\tau}{n}} = c_\varepsilon.
\]

\n It follows that

\[
J_\varepsilon(\tilde{u}_\varepsilon) = c_\varepsilon \; .
\]

Now, we will show that

\[
\int_M |\nabla_g \tilde{u}_\varepsilon|^p dv_g \not = 0 \; .
\]

\n Since $\nabla_g \tilde{u}_\varepsilon = \pm \nabla_g
\tilde{u}_\varepsilon$ almost everywhere and towards a future
simplification of notation we can assume $\tilde{u}_\varepsilon
\geq 0$. Suppose by contradiction that $\tilde{u}_\varepsilon$ is
constant. As

\[
\int_M |\tilde{u}_\varepsilon|^r dv_g = 1 \; ,
\]

\n we must have $\tilde{u}_\varepsilon = |M|^{- \frac{1}{r}}$. Also we
have

\[
1 > c_\varepsilon = \left( A(p,n)^{\frac{\tau}{p}} \left(\int_M |\nabla_g
\tilde{u}_\varepsilon|^p\; dv_g\right)^{\frac{\tau}{p}} + (B_\varepsilon -
\gamma_\varepsilon) \left(\int_M |\tilde{u}_\varepsilon|^p\;
dv_g\right)^{\frac{\tau}{p}} \right) \left( \int_M
|\tilde{u}_\varepsilon|^p\; dv_g \right)^{\frac{\tau}{n}} \; .
\]

\n As $\tilde{u}_\varepsilon = |M|^{-\frac{1}{r}}$, we obtain

\[
1 > (B_\varepsilon - \gamma_\varepsilon) |M|^{\frac{\tau}{n}} \; .
\]

\n By our choice of $\gamma_\varepsilon$ we have that $B_\varepsilon -
\gamma_\varepsilon > |M|^{- \frac{\tau}{n}}$ what produces the
contradiction. Therefore $\tilde{u}_\varepsilon$ is not
constant.

\n Define now

\[
v_\varepsilon = \frac{\tilde{u}_\varepsilon}{||\nabla_g
\tilde{u}_\varepsilon||_{L^p(M)}} \; .
\]

\n Since $\tilde{u}_\varepsilon \geq 0$, we have
$v_\varepsilon \geq 0$ in $M$. Note that $v_\varepsilon$ satisfies

\[
\left(A(p,n)^{\frac{\tau}{p}} + (B_\varepsilon - \gamma_\varepsilon) \left(\int_M
|v_\varepsilon|^p dv_g \right)^{\frac{\tau}{p}}\right)
\left(\int_M |v_\varepsilon|^p dv_g \right)^{\frac{\tau}{n}} =
c_\varepsilon \left(\int_M |v_\varepsilon|^r dv_g
\right)^{\frac{\tau}{p}}
\]

\n or

\[
\frac{A(p,n)^{\frac{\tau}{p}}}{c_\varepsilon} = \left(\int_M|v_\varepsilon|^r
dv_g\right)^{\frac{\tau}{p}} \left( \int_M |v_\varepsilon|^p dv_g
\right)^{- \frac{\tau}{n}} - \frac{B_\varepsilon -
\gamma_\varepsilon}{c_\varepsilon} \left(\int_M |v_\varepsilon|^p
dv_g \right)^{\frac{\tau}{p}} .
\]

\n At last, define

\[
G_\varepsilon(u) = \left(\int_M|u|^r dv_g\right)^{\frac{\tau}{p}}
\left( \int_M |u|^p dv_g \right)^{- \frac{\tau}{n}} -
\frac{B_\varepsilon - \gamma_\varepsilon}{c_\varepsilon} \left(\int_M |u|^p
dv_g \right)^{\frac{\tau}{p}} .
\]

\n By (\ref{um}), we have that $G_\varepsilon(u) \leq
\frac{A(p,n)^{\frac{\tau}{p}}}{c_\varepsilon}$ for all functions $u
\in {\cal D} = \{u\in H^{1,p}(M); \int_M |\nabla_g u|^p dv_g =
1\}$. Since $G_\varepsilon$ is of class $C^1$, $v_\varepsilon
\in {\cal D}$, $c_\varepsilon < 1$ and
$G_\varepsilon(v_\varepsilon) =
\frac{A(p,n)^{\frac{\tau}{p}}}{c_\varepsilon}$, we have

\begin{equation}\label{3dha}
\nu_\varepsilon = \sup_{u \in {\cal D}} G_\varepsilon (u) =
G_\varepsilon(v_\varepsilon) > A(p,n)^{\frac{\tau}{p}}\, .
\end{equation}

\n So, $v_\varepsilon$ satisfies the Euler-Lagrange
equation

\begin{equation}\label{eq1}
\frac{n + p}{n}||v_\varepsilon||_{L^r(M)}^{\frac{\tau(n +
p) - r n}{n}} ||v_\varepsilon||_{L^p(M)}^{- \frac{\tau
p}{n}} v_\varepsilon^{r - 1} - \frac{p}{n}
||v_\varepsilon||_{L^r(M)}^{\frac{\tau(n + p)}{n}}
||v_\varepsilon||_{L^p(M)}^{-\frac{\tau p}{n} - p}
v_\varepsilon^{p - 1} - d_\varepsilon
||v_\varepsilon||_{L^p(M)}^{\tau - p}
v_\varepsilon^{p - 1} = \nu_\varepsilon \Delta_{p,g}
v_\varepsilon \; ,
\end{equation}

\n where $\Delta_{p,g} = -{\rm div}_g(|\nabla_g|^{p-2} \nabla_g)$
is the $p$-Laplace operator of $g$ and $d_\varepsilon =
\frac{B_\varepsilon - \gamma_\varepsilon}{c_\varepsilon}$.

We now set

\[
u_\varepsilon = \frac{v_\varepsilon}{||v_\varepsilon||_{L^r(M)}} \geq 0 \;.
\]

\n Writing the Euler-Lagrange equation in terms of $u_\varepsilon$,
we have

\begin{equation} \label{3ep}
\lambda_\varepsilon^{-1} A_\varepsilon \Delta_{p,g} u_\varepsilon
+ d_\varepsilon
 A_\varepsilon^{\frac{\tau}{p}} ||u_\varepsilon||_{L^p(M)}^{\tau -
p} u_\varepsilon^{p - 1} + \frac{p}{n}
||u_\varepsilon||_{L^p(M)}^{-p} u_\varepsilon^{p - 1} = \frac{n +
p}{n} u_\varepsilon^{r - 1}\ \ {\rm on}\ \ M\, ,
\end{equation}

\n where  $||u_\varepsilon||_{L^r(M)} = 1$,

\[
A_\varepsilon = \left(\int_M u_\varepsilon^p\; dv_g
\right)^{\frac{p}{n}}
\]

\n and

\[
\lambda_\varepsilon = \nu_\varepsilon^{-1}
||u_\varepsilon||_{L^p(M)}^{\frac{p(p - \tau)}{n}}
||v_\varepsilon||_{L^r(M)}^{\tau - p} \; .
\]

\n We now highlight an important consequence of the Euler-Lagrange
equation:

\begin{equation}\label{rm}
\lambda_\varepsilon^{-1} \geq A(p,n) \; .
\end{equation}

\n In order to show (\ref{rm}), note that taking $v_\varepsilon$ as the
test function in (\ref{eq1}), we have

\[
\nu_\varepsilon \leq ||v_\varepsilon||_{L^r(M)}^{\frac{\tau(n +
p)}{n}} ||v_\varepsilon||_{L^p(M)}^{- \frac{\tau p}{n}} \; .
\]

\n Putting together the previous inequality, (\ref{3dha}) and noting that $\tau \leq p$, we get

\[
\lambda_\varepsilon = \nu_\varepsilon^{-1}
||v_\varepsilon||_{L^r(M)}^{\frac{p(\tau - p)}{n}}
||v_\varepsilon||_{L^p(M)}^{-\frac{p(\tau - p)}{n}}
||v_\varepsilon||_{L^r(M)}^{\tau - p} = \nu_\varepsilon^{-1}
\left(||v_\varepsilon||_{L^r(M)}^{\frac{\tau(n + p)}{n}}
||v_\varepsilon||_{L^p(M)}^{- \frac{\tau p}{n}}
\right)^{\frac{\tau - p}{\tau}} \leq
\nu_\varepsilon^{-\frac{p}{\tau}} \leq A(p,n)^{-1}
\; .
\]

\n This proves (\ref{rm}).

As $p < r$, there exists a constant $c > 0$ such that

\[
0 \leq A_\varepsilon \leq c \; ,
\]

\n for all $\varepsilon > 0$.

\n Now, we will examine two possible situations:

{\bf A:}

\[
\lim_{\varepsilon \rightarrow 0} A_\varepsilon > 0
\]

{\bf B:}

\begin{equation}\label{contradicao}
\lim_{\varepsilon \rightarrow 0} A_\varepsilon = 0 \; .
\end{equation}

\n If occur (A), we have by (\ref{3ep}) that $d_\varepsilon$ is
bounded, consequently, ${\cal B}$ is bounded.

The remainder of the proof will be to show that the only
possibility of (B) to occur is when $\tau = \min\{2,p\}$, and in
that case, we should have also ${\cal B}$ is limited. The argument in
this case is more delicate, so we divide the proof in several
steps. The central idea is around a study of explosion.

\n Assume that (B) occurs. First we observe, using (\ref{epsilon})
with $\varepsilon = \varepsilon_0$ and because $p < r$, that there exists
a constant $c > 0$ such that

\[
A_\varepsilon \int_M |\nabla_g u_\varepsilon|^p dv_g
> c \;.
\]

\n Using this in (\ref{3ep}), we have that

\begin{equation}\label{lambda}
\lambda_\epsilon^{-1} \leq c
\end{equation}

\n for some $c > 0$.

\subsection{$L^r$-concentration} \label{concentration}

Using (\ref{3ep}) and by Tolksdorf's regularity theory (see \cite{To}), it follows
that $u_\varepsilon$ is of class $C^1$. So, let $x_\varepsilon \in M$ be a maximum point of $u_\varepsilon$,
that is,

\begin{equation}\label{3ix}
u_\varepsilon(x_\varepsilon) = ||u_\varepsilon||_{L^\infty(M)}\, .
\end{equation}

\n Our aim here is to establish that

\begin{equation} \label{4}
\lim_{\sigma \rightarrow \infty}\lim_{\varepsilon \rightarrow 0}
\int_{B(x_\varepsilon,\sigma A_\varepsilon^{\frac{n + p}{p^2}})}
u_\varepsilon^r\; dv_g = 1\, ,
\end{equation}

\n Note that $A_\varepsilon^{\frac{n + p}{p^2}} \rightarrow 0$ as
$\varepsilon \rightarrow 0$.

\n For each $x \in B(0, \sigma)$, define

\begin{equation}\label{l}
\begin{array}{l}
h_\varepsilon(x) = g(\exp_{x_\varepsilon} (A_\varepsilon^{\frac{n + p}{p^2}} x))\, , \vspace{0,3cm}\\
\varphi_\varepsilon(x) = A_\varepsilon^{\frac{n^2}{p^3}}
u_\varepsilon(\exp_{x_\varepsilon} (A_\varepsilon^{\frac{n +
p}{p^2}} x))\, .
\end{array}
\end{equation}

\n By (\ref{3ep}), one easily deduces that

\begin{equation}\label{eqgn}
\lambda_\varepsilon^{-1} \Delta_{p,h_\varepsilon}
\varphi_\varepsilon + d_\varepsilon A_\varepsilon^{\frac{\tau}{p}
- 1} ||u_\varepsilon||_{L^p(M)}^{\tau - p} A_\varepsilon^{\frac{n
+ p}{p}} \varphi_\varepsilon^{p - 1} + \frac{p}{n}
\varphi_\varepsilon^{p - 1} = \frac{n + p}{n}
\varphi_\varepsilon^{r - 1}\ \ {\rm on}\ \ B(0, \sigma)\, .
\end{equation}

\n For $p < n$, using (\ref{rm}) and applying the Moser's
iterative scheme (see \cite{serrin}) to this last equation, we see
that

\[
A_\varepsilon^{\frac{n(n + p)}{p^2}}
||u_\varepsilon||_{L^\infty(M)}^r = \sup_{B(0,\frac{\sigma}{2})}
\varphi_\varepsilon^r \leq c \int_{B(0,\sigma)}
\varphi_\varepsilon^r\; dh_\varepsilon = c \int_{B(x_\varepsilon,
\sigma A_\varepsilon^{\frac{n + p}{p^2}})} u_\varepsilon^r\; dv_g
\leq c \; ,
\]

\n for $\varepsilon$ small enough. For $p \geq n$, we use
(\ref{rm}) and Morrey's inequality (see \cite{serrin} for general
idea) also for

\[
A_\varepsilon^{\frac{n(n + p)}{p^2}} ||u_\varepsilon||_{L^\infty(M)}^r = \sup_{B(0,\frac{\sigma}{2})} \varphi_\varepsilon^r \leq c \; .
\]

\n This estimates together with

\[
1 = \int_M u_\varepsilon^r\; dv_g \leq
||u_\varepsilon||_{L^\infty(M)}^{r-p} \int_M u_\varepsilon^p dv_g
= \left(||u_\varepsilon||_{L^\infty(M)} \;
A_\varepsilon^{\frac{n^2}{p^3}}\right)^{r-p}
\]

\n gives

\begin{equation}\label{3dqp}
1 \leq ||u_\varepsilon||_{L^\infty(M)} \;
A_\varepsilon^{\frac{n^2}{p^3}} \leq c \, ,
\end{equation}

\n for $\varepsilon > 0$ and $p > 1$.

\n We will use in the sequel that (\ref{3dqp}) means that the
limiting behavior of $||u_\varepsilon||_{L^\infty(M)} $ and $
A_\varepsilon^{-\frac{n^2}{p^3}}$ have the same order.

In particular, there exists a constant $c > 0$ such that

\begin{equation}\label{lb}
\int_{B(0,\sigma)} \varphi_\varepsilon^r\; dh_\varepsilon \geq c >
0
\end{equation}

\n for $\varepsilon$ small enough.

\n Now using Cartan expansion in normal coordinates and (\ref{3dqp}), we have for each $\sigma > 0$, that

\[
\int_{B(0,\sigma)} \varphi_\varepsilon^p dx \leq c
\int_{B(0,\sigma)} \varphi_\varepsilon^p dh_\varepsilon = c \;
\frac{\int_{B(x_\varepsilon,\sigma A_\varepsilon^{\frac{n +
p}{p^2}})} u_\varepsilon^p\; dv_g}{\int_M u_\varepsilon^p\; dv_g}
\leq c \; ,
\]

\n where $c$ is independent of $\varepsilon$. We have also

\[
\int_{B(0,\sigma)} |\nabla \varphi_\varepsilon|^p\; dx \leq c
\int_{B(0,\sigma)} |\nabla_{h_\varepsilon}
\varphi_\varepsilon|^p\; dh_\varepsilon = c A_\varepsilon
\int_{B(x_\varepsilon,\sigma A_\varepsilon^{\frac{n + p}{p^2}})}
|\nabla_g u_\varepsilon|^p\; dv_g \leq c \; A(p,n)^{-1}\, .
\]

\n Therefore there exists $\varphi \in W^{1,p}(\R^n)$ such that,
for some subsequence, $\varphi_\varepsilon \rightharpoonup
\varphi$ in $W^{1,p}_{loc}(\R^n)$. For each $\sigma > 0$, we have

\[
\int_{B(0,\sigma)} \varphi^r\; dx = \lim_{\varepsilon \rightarrow
0} \int_{B(0,\sigma)} \varphi_\varepsilon^r\; dh_\varepsilon =
\lim_{\varepsilon \rightarrow 0} \int_{B(x_\varepsilon,\sigma
A_\varepsilon^{\frac{n + p}{p^2}})} u_\varepsilon^r\; dv_g \leq 1
\; .
\]

\n In particular,

\begin{equation}\label{inter}
\varphi \in L^r(\R^n) \; .
\end{equation}

\n Let $\eta \in C_0^1(\R)$ be a cutoff function such that $\eta =
1$ on $[0,\frac{1}{2}]$, $\eta = 0$ on $[1,\infty)$ and $0 \leq
\eta \leq 1$. Define $\eta_{\varepsilon,\sigma}(x) = \eta((\sigma
A_\varepsilon^{\frac{n + p}{p^2}})^{-1} d_g(x,x_\varepsilon))$.
Choosing $u_\varepsilon \eta_{\varepsilon,\sigma}^p$ as a test
function in (\ref{3ep}), one gets

\[
\lambda_\varepsilon^{-1} A_\varepsilon \int_M |\nabla_g
u_\varepsilon|^p \eta_{\varepsilon,\sigma}^p\; dv_g +
\lambda_\varepsilon^{-1} A_\varepsilon \int_M |\nabla_g
u_\varepsilon|^{p - 2} \nabla_g u_\varepsilon \cdot \nabla_g
(\eta_{\varepsilon,\sigma}^p) u_\varepsilon\; dv_g + \frac{p}{n}
\frac{\int_M u_\varepsilon^p \eta_{\varepsilon,\sigma}^p\;
dv_g}{\int_M u_\varepsilon^p dv_g}
\]

\begin{equation}\label{3sc3}
\leq \frac{n + p}{n} \int_M u_\varepsilon^r
\eta_{\varepsilon,\sigma}^p\; dv_g\, .
\end{equation}

\n We now show that

\begin{equation}\label{ef1}
\lim_{\sigma \rightarrow \infty} \lim_{\varepsilon \rightarrow 0} A_\varepsilon \int_M
|\nabla_g u_\varepsilon|^{p - 2} \nabla_g u_\varepsilon \cdot
\nabla_g (\eta_{\varepsilon,\sigma}^p) u_\varepsilon\; dv_g = 0\,
.
\end{equation}

\n Taking $u_\varepsilon$ as the test function, by (\ref{3ep}) we have

\[
A_\varepsilon \int_M |\nabla_g u_\varepsilon|^p dv_g \leq
\lambda_\varepsilon \leq A(p,n)^{-1} \; .
\]

\n Therefore, it suffices to establish that

\begin{equation}\label{ef6}
A_\varepsilon \int_M u_\varepsilon^p |\nabla_g
\eta_{\varepsilon,\sigma}|^p\; dv_g \leq \frac{c}{\sigma^p}\, .
\end{equation}

\n Using (\ref{3dqp}), we derive

\[
A_\varepsilon \int_M u_\varepsilon^p |\nabla_g
\eta_{\varepsilon,\sigma}|^p\; dv_g \leq c \frac{A_\varepsilon
||u_\varepsilon||_{L^\infty(M)}^{p - r}}{\sigma^p
A_\varepsilon^{\frac{n + p}{p}}} \int_M u_\varepsilon^r\; dv_g
\leq c \frac{A_\varepsilon A_\varepsilon^{\frac{n}{p}}}{\sigma ^p
A_\varepsilon^{\frac{n + p}{p}}} \leq \frac{c}{\sigma^p}.
\]

\n Therefore (\ref{ef6}) holds and (\ref{ef1}) is valid.

\n Replacing (\ref{rm}) and (\ref{ef1}) in (\ref{3sc3}), one
arrives at

\[
\frac{n}{n + p} \lim_{\sigma \rightarrow
\infty} \lim_{\varepsilon \rightarrow 0}\left(A(p,n)\; A_\varepsilon \int_M
|\nabla_g u_\varepsilon|^p \eta_{\varepsilon, \sigma}^p\; dv_g
\right) + \frac{p}{n + p} \lim_{\sigma \rightarrow
\infty}\lim_{\varepsilon \rightarrow 0} \frac{\int_M u_\varepsilon^p \eta_{\varepsilon,\sigma}^p\;
dv_g}{\int_M u_\varepsilon^p\; dv_g} \leq \lim_{\sigma,
 \rightarrow \infty} \lim_{\varepsilon \rightarrow 0} \int_M u_\varepsilon^r
\eta_{\varepsilon, \sigma}^p\; dv_g\, .
\]

\n To rewrite this inequality in a suitable format, we remark that

\[
\left| \int_M u_\varepsilon^r \eta_{\varepsilon,\sigma}^p \; dv_g
- \int_M u_\varepsilon^r \eta_{\varepsilon,\sigma}^r \;dv_g
\right| \leq \int_{B(x_\varepsilon,\sigma A_\varepsilon^{\frac{n +
p}{p^2}}) \setminus B(x_\varepsilon, (\sigma
A_\varepsilon^{\frac{n + p}{p^2}})/2)} u_\varepsilon^r \;dv_g =
\int_{B(0,\sigma) \setminus B(0,\sigma/2)} \varphi_\varepsilon^r
dh_\varepsilon \; ,
\]

\n using (\ref{inter}), so that

\[
\lim_{\sigma \rightarrow \infty} \lim_{\varepsilon \rightarrow 0}\int_M
u_\varepsilon^r \eta_{\varepsilon,\sigma}^r\; dv_g =
\lim_{\sigma \rightarrow \infty} \lim_{\varepsilon \rightarrow 0} \int_M
u_\varepsilon^r \eta_{\varepsilon,\sigma}^p\; dv_g\, .
\]

\n Consequently, we can write

\begin{equation}\label{ef2}
\frac{n}{n + p} \lim_{\sigma \rightarrow
\infty} \lim_{\varepsilon \rightarrow 0}\left(A(p,n) \; A_\varepsilon \int_M
|\nabla_g u_\varepsilon|^p \eta_{\varepsilon, \sigma}^p\; dv_g
\right) + \frac{p}{n + p} \lim_{\sigma \rightarrow
\infty} \lim_{\varepsilon \rightarrow 0} \frac{\int_M u_\varepsilon^p \eta_{\varepsilon,\sigma}^p\;
dv_g}{\int_M u_\varepsilon^p\; dv_g} \leq \lim_{\sigma \rightarrow \infty} \lim_{\varepsilon \rightarrow 0} \int_M u_\varepsilon^r
\eta_{\varepsilon, \sigma}^r\; dv_g\, .
\end{equation}

\n On the other hand, for $\kappa > 0$ let the constant
$B_\kappa > 0$, independent of $\varepsilon$, such that

\[
\left(\int_M u_\varepsilon^r \eta^r_{\varepsilon,\sigma}\; dv_g
\right)^{\frac{\tau}{p}} \leq \left( (A(p,n)^{\frac{\tau}{p}} + \kappa)
\left(\int_M |\nabla_g(u_\varepsilon
\eta_{\varepsilon,\sigma})|^p\; dv_g\right)^{\frac{\tau}{p}} +
B_\kappa \left(\int_M u_\varepsilon^p
\eta_{\varepsilon,\sigma}^p\; dv_g\right)^{\frac{\tau}{p}} \right)
\left(\int_M u_\varepsilon^p \eta^p_{\varepsilon,\sigma}\; dv_g
\right)^{\frac{\tau}{n}}\, .
\]

\n From the definition of $A_\varepsilon$, Young's inequality and
$(x + y)^p \leq x^p + c x^{p - 1}y + c y^p$ for $x,y \geq 0$, one
has

\[
\left(\int_M u_\varepsilon^r \eta^r_{\varepsilon,\sigma} dv_g\right)^{\frac{\tau}{p}}
\leq c \left(A_\varepsilon \int_M
u_\varepsilon^p dv_g \right)^{\frac{\tau}{p}}
\]
\[
+ (A(p,n)^{\frac{\tau}{p}} + \kappa) \left(\left((1 + \kappa) \int_M
|\nabla_g u_\varepsilon|^p \eta_{\varepsilon,\sigma}^p\; dv_g +
c(\kappa) \int_M u_\varepsilon^p |\nabla_g
\eta_{\varepsilon,\sigma}|^p dv_g \right)^{\frac{\tau}{p}}\right)
\left(\int_M u_\varepsilon^p \eta_{\varepsilon,\sigma}^p\; dv_g
\right)^{\frac{\tau}{n}} \, .
\]

\n Then, using hypothesis (B) and (\ref{ef6}), letting
$\varepsilon \rightarrow 0$, $\sigma \rightarrow \infty$ and
$\kappa \rightarrow 0$, one gets

\begin{equation}\label{3il2}
\lim_{\sigma \rightarrow \infty} \lim_{\varepsilon \rightarrow 0}  \int_M
u_\varepsilon^r \eta_{\varepsilon,\sigma}^r\; dv_g
\leq \lim_{\sigma \rightarrow \infty} \lim_{\varepsilon \rightarrow 0}\left( A(p,n) A_\varepsilon
\int_M |\nabla_g u_\varepsilon|^p \eta_{\varepsilon,\sigma}^p\;
dv_g \right) \lim_{\sigma \rightarrow \infty} \lim_{\varepsilon \rightarrow 0}
\left(\frac{\int_M u_\varepsilon^p \eta_{\varepsilon,\sigma}^p\;
dv_g}{\int_M u_\varepsilon^p\; dv_g} \right)^{\frac{p}{n}}\, .
\end{equation}

\n Let

\[
X = \lim_{\sigma \rightarrow \infty} \lim_{\varepsilon \rightarrow 0} \left(
A(p,n) \; A_\varepsilon \; \int_M |\nabla_g
u_\varepsilon|^p \eta_{\varepsilon,\sigma}^p dv_g \right), \ \ Y =
\lim_{\sigma \rightarrow \infty} \lim_{\varepsilon \rightarrow 0}\frac{\int_M
u_\varepsilon^p \eta_{\varepsilon,\sigma}^p dv_g}{\int_M
u_\varepsilon^p dv_g}, \ \ Z = \lim_{\sigma \rightarrow \infty} \lim_{\varepsilon \rightarrow 0} \int_M u_\varepsilon^r
\eta_{\varepsilon,\sigma}^r dv_g\, .
\]

\n Clearly, $X, Y, Z \leq 1$ and (\ref{ef2}) and (\ref{3il2}) may be rewritten as

\begin{equation}\label{3sc4}
\left\{%
\begin{array}{c}
\frac{n}{n + p} X + \frac{p}{n + p} Y \leq Z   \vspace{0,5cm}\\
Z \leq X Y^{\frac{p}{n}} \, .
\end{array}%
\right.
\end{equation}

\n By (\ref{lb}), one also has $Z > 0$, so that $X, Y > 0$.

 We  are now ready to prove that $Z = 1$. By Young's inequality, we have that

\[
X^{\frac{n}{n + p}} Y^{\frac{p}{n + p}} \leq Z
\]

\n so

\[
1 \leq X^{\frac{p}{n + p}} Y^{\frac{p^2}{n(n + p)}}.
\]

This implies $X = Y =1$, therefore $Z=1$. This prove (\ref{4}).

\subsection{Pointwise estimates} \label{pw}

There exists a constant $c > 0$, independent of $\varepsilon$, such that

\[
d_g(x,x_\epsilon) u_\epsilon(x)^{\frac{p}{n}} \leq c
A_\epsilon^{\frac{1}{p}}
\]

\n \emph{for all $x \in M$ and $\epsilon > 0$ small enough.}

\n The proof proceeds by contradiction. Suppose that the assertion
above is false. Then, there exists $y_\varepsilon \in M$ such that
$f_\varepsilon(y_\varepsilon) \rightarrow \infty$ as $\varepsilon
\rightarrow 0$, where

\[
f_\varepsilon(x) = d_g(x_\varepsilon,x)
u_\varepsilon(x)^{\frac{p}{n}} A_\varepsilon^{-\frac{1}{p}} \; .
\]

\n Assume, without loss of generality, that
$f_{\varepsilon}(y_\varepsilon) =
||f_{\varepsilon}||_{L^\infty(M)}$. From (\ref{3dqp}), we have

\[
f_\varepsilon(y_\varepsilon) \leq
\left(\frac{u_\varepsilon(y_\varepsilon)}{||u_\varepsilon||_{L^\infty(M)}}\right)^{\frac{p}{n}}
d_g(x_\varepsilon,y_\varepsilon)
||u_\varepsilon||_{L^\infty(M)}^{\frac{p(n + p)}{n^2}} \leq c
d_g(x_\varepsilon,y_\varepsilon
||u_\varepsilon||_{L^\infty(M)}^{\frac{p(n + p)}{n^2}} \, ,
\]

\n so that
\begin{equation}\label{ldi}
d_g(x_\varepsilon,y_\varepsilon)
||u_\varepsilon||_{L^\infty(M)}^{\frac{p(n + p)}{n^2}} \rightarrow
\infty\, .
\end{equation}

\n For any $\sigma > 0$ fixed and $\varepsilon \in (0,1)$, we show that

\begin{equation}\label{int}
B(y_\varepsilon,\varepsilon d_g(x_\varepsilon,y_\varepsilon)) \cap
B(x_\varepsilon, \sigma
||u_\varepsilon||_{L^\infty(M)}^{-\frac{p(n + p)}{n^2}}) = \emptyset
\end{equation}

\n for $\varepsilon > 0$ small enough. Clearly, this assertion
follows from

\[
d_g(x_\varepsilon,y_\varepsilon) \geq \sigma
||u_\varepsilon||_{L^\infty(M)}^{-\frac{p(n + p)}{n^2}} + \varepsilon
d_g(x_\varepsilon,y_\varepsilon) \, .
\]

\n But the above inequality is equivalent to
\[
(1 - \varepsilon)d_g(x_\varepsilon, y_\varepsilon)
||u_\varepsilon||_{L^\infty(M)}^{\frac{p(n + p)}{n^2}} \geq \sigma \, ,
\]

\n which is clearly satisfied, since
$d_g(x_\varepsilon,y_\varepsilon)
||u_\varepsilon||_{L^\infty(M)}^{\frac{p(n + p)}{n^2}} \rightarrow
\infty$ and $1 - \varepsilon > 0$. Then it follows (\ref{int}).

\n We claim that there exists a constant $c > 0$ such that

\begin{equation}\label{lim}
u_\varepsilon(x) \leq c u_\varepsilon(y_\varepsilon)
\end{equation}

\n for all $x \in B(y_\varepsilon, \varepsilon
d_g(x_\varepsilon,y_\varepsilon))$ and $\varepsilon > 0$ small
enough. Indeed, for each $x \in B(y_\varepsilon, \varepsilon
d_g(x_\varepsilon,y_\varepsilon))$, we have

\[
d_g(x,x_\varepsilon) \geq d_g(x_\varepsilon,y_\varepsilon) -
d_g(x,y_\varepsilon) \geq (1 - \varepsilon)
d_g(x_\varepsilon,y_\varepsilon)\, .
\]

\n Thus,

\[
d_g(y_\varepsilon,x_\varepsilon)
u_\varepsilon(y_\varepsilon)^{\frac{p}{n}}
A_\varepsilon^{-\frac{1}{p}} = f_{\varepsilon}(y_\varepsilon) \geq
f_{\varepsilon}(x) = d_g(x,x_\varepsilon)
u_\varepsilon(x)^{\frac{p}{n}} A_\varepsilon^{-\frac{1}{p}}
\]

\[
\geq (1 - \varepsilon) d_g(y_\varepsilon,x_\varepsilon)
u_\varepsilon(x)^{\frac{p}{n}} A_\varepsilon^{-\frac{1}{p}}\, ,
\]

\n so that

\[
u_\varepsilon(x) \leq \left(\frac{1}{1 -
\varepsilon}\right)^{\frac{n}{p}} u_\varepsilon(y_\varepsilon)
\]

\n for all $x \in B(y_\varepsilon, \varepsilon
d_g(x_\varepsilon,y_\varepsilon))$ and $\varepsilon > 0$ small
enough. This proves our claim.

\n Since

\[
d_g(x_\varepsilon,y_\varepsilon)
u_\varepsilon(y_\varepsilon)^{\frac{p}{n}}
A_\varepsilon^{-\frac{1}{p}} \rightarrow \infty
\]

\n we have

\[
u_\varepsilon(y_\varepsilon)^{- \frac{p}{n}}
A_\varepsilon^{\frac{1}{p}} \rightarrow 0,
\]

\n because $M$ is compact.

\n So, we can define

\[
\begin{array}{l}
h_\varepsilon(x) =
g(\exp_{y_\varepsilon}(A_\varepsilon^{\frac{1}{p}}
u_\varepsilon(y_\varepsilon)^{-\frac{p}{n}} x)), \vspace{0,2cm}\\
\psi_\varepsilon(x) = u_\varepsilon(y_\varepsilon)^{-1}
u_\varepsilon(\exp_{y_\varepsilon}(A_\varepsilon^{\frac{1}{p}}
u_\varepsilon(y_\varepsilon)^{-\frac{p}{n}} x))
\end{array}
\]

\n for each $x \in B(0,2)$ and $\varepsilon$ small enough.

\n From (\ref{3ep}), it follows that

\begin{equation}\label{3sc5}
\lambda_\varepsilon^{-1} \Delta_{p,h_\varepsilon} \psi_\varepsilon
+ d_\varepsilon A_\varepsilon^{\frac{\tau}{p}}
||u_\varepsilon||_{L^p(M)}^{\tau - p}
u_\varepsilon(y_\varepsilon)^{p - r} \psi_\varepsilon^{p - 1} +
\frac{p}{n} ||u_\varepsilon||_{L^p(M)}^{-p}
u_\varepsilon(y_\varepsilon)^{p - r} \psi_\varepsilon^{p - 1} =
\frac{n + p}{n} \psi_\varepsilon^{r - 1} \ \ {\rm on}\ \ B(0,2)\,
.
\end{equation}

\n In particular,

\[
\int_{B(0,2)} |\nabla_{h_\varepsilon} \psi_\varepsilon|^{p - 2} \nabla_{h_\varepsilon} \psi_{\varepsilon} \cdot\nabla_{h_\varepsilon}\phi\; dv_{h_\varepsilon} \leq c \int_{B(0,2)} \psi_\varepsilon^{r - 1} \phi\; dv_{h_\varepsilon}
\]

\n for all positive test function $\phi \in C_0^1(B(0,2))$. So, by
(\ref{3dqp}) and the Moser's iterative scheme if $p < n$ or
Morrey's inequality if $p \geq n$ (see \cite{serrin}), one deduces
that

\[
1 \leq \sup_{B(0,\frac{1}{4})} \psi_\varepsilon^r \leq c
\int_{B(0,\frac{1}{2})} \psi_\varepsilon^r\; dv_{h_\varepsilon} =
c \left( A_\varepsilon^{\frac{n}{p}}
u_\varepsilon(y_\varepsilon)^{\frac{p^2}{n}} \right)^{- 1}
\int_{B(y_\varepsilon, \frac{1}{2} A_\varepsilon^{\frac{1}{p}}
u_\varepsilon(y_\varepsilon)^{-\frac{p}{n}})} u_\varepsilon^r\;
dv_g
\]
\[
\leq c
\left(\frac{||u_\varepsilon||_{L^\infty(M)}}{u_\varepsilon(y_\varepsilon)}
\right)^{\frac{p^2}{n}} \int_{B(y_\varepsilon,
\frac{1}{2} A_\varepsilon^{\frac{1}{p}}
u_\varepsilon(y_\varepsilon)^{-\frac{p}{n}})} u_\varepsilon^r\;
dv_g\,.
\]

\n For simplicity, we rewrite this last inequality as

\begin{equation}\label{3dcm}
0 < c \leq m_\varepsilon^{\frac{p^2}{n}} \int_{B(y_\varepsilon,\frac{1}{2}
A_\varepsilon^{\frac{1}{p}} u_\varepsilon(y_\varepsilon)^{-\frac{p}{n}})} u_\varepsilon^r\;
dv_g\, ,
\end{equation}

\n where $m_\varepsilon =
\frac{||u_\varepsilon||_{L^\infty(M)}}{u_\varepsilon(y_\varepsilon)}$.

\n Observe that $B(y_\varepsilon, \frac{1}{2}
A_\varepsilon^{\frac{1}{p}}
u_\varepsilon(y_\varepsilon)^{-\frac{p}{n}}) \subset
B(y_\varepsilon,\varepsilon d(x_\varepsilon,y_\varepsilon))$ for
$\varepsilon$ small enough, this follows directly from the
contradiction hypothesis. So, the $L^r$-concentration property
(\ref{4}) combined with (\ref{int}) provides

\[
\int_{B(y_\varepsilon, \frac{1}{2} A_\varepsilon^{\frac{1}{p}}
u_\varepsilon(y_\varepsilon)^{-\frac{p}{n}})} u_\varepsilon^r\;
dv_g \rightarrow 0 \; ,
\]

\n when $\varepsilon \rightarrow 0$.

\n So, we have that

\[
\lim_{\varepsilon \rightarrow 0} m_\varepsilon = \infty \;.
\]

\n Our goal now is to contradict (\ref{3dcm}).
Initially, from (\ref{3dqp}) and (\ref{lim}), we have

\begin{equation}\label{3dpt}
m_\varepsilon^{\frac{p^2}{n}} \int_{D_{\varepsilon}}
u_\varepsilon^r\; dv_g \leq c \; m_\varepsilon^{\frac{p^2}{n}}
||u_\varepsilon||^r_{L^\infty(D_{\varepsilon})}
(A_\varepsilon^{\frac{1}{p}}
u_\varepsilon(y_\varepsilon)^{-\frac{p}{n}})^n \leq c \;
m_\varepsilon^{\frac{p^2}{n}} u_\varepsilon(y_\varepsilon)^r
(A_\varepsilon^{\frac{1}{p}}
u_\varepsilon(y_\varepsilon)^{-\frac{p}{n}})^n \leq c\, ,
\end{equation}

\n where $D_{\varepsilon} = B(y_\varepsilon,
A_\varepsilon^{\frac{1}{p}}
u_\varepsilon(y_\varepsilon)^{-\frac{p}{n}})$ \; .

\n Consider the function $\eta_\varepsilon(x) = \eta(
A_\varepsilon^{-\frac{1}{p}} d_g(x,y_\varepsilon)
u_\varepsilon(y_\varepsilon)^{\frac{p}{n}})$, where $\eta \in
C_0^1(\R)$ is a cutoff function such that $\eta = 1$ on
$[0,\frac{1}{2}]$, $\eta = 0 $ on $[1,\infty)$ and $0 \leq \eta
\leq 1$. Taking $u_\varepsilon \eta_\varepsilon^p$ as a test
function in (\ref{3ep}), one has

\[
\lambda_\varepsilon ^{-1} A_\varepsilon \int_M |\nabla_g
u_\varepsilon|^p \eta_\varepsilon^p\; dv_g +
\lambda_\varepsilon^{-1} p A_\varepsilon \int_M |\nabla_g
u_\varepsilon|^{p - 2} u_{\varepsilon}\eta_{\varepsilon}^{p-1} \nabla_g u_\varepsilon\cdot\nabla_g\eta_{\varepsilon}
\; dv_g + d_\varepsilon
A_\varepsilon^{\frac{\tau}{p}} ||u_\varepsilon||_{L^p(M)}^{\tau -
p} \int_M u_\varepsilon^p \eta_\varepsilon^p\; dv_g
\]
\[
+ \frac{p}{n} \frac{\int_M u_\varepsilon^p
\eta_\varepsilon^p}{\int_M u_\varepsilon^p dv_g}\; dv_g =
\frac{n + p}{n} \int_M u_\varepsilon^r \eta_\varepsilon^p\;
dv_g\, .
\]

\n By H\"{o}lder and Young's inequalities,
\[
\left| \int_M |\nabla_g u_\varepsilon|^{p - 1}  \nabla_g \eta_\varepsilon u_\varepsilon
\eta_\varepsilon^{p - 1} \; dv_g \right| \leq \delta \int_M
|\nabla_g u_\varepsilon|^p \eta_\varepsilon^p\; dv_g + c(\delta)
\int_M |\nabla_g \eta_\varepsilon|^p u_\varepsilon^p\; dv_g\, .
\]

\n Also, by (\ref{3dqp}) and (\ref{lim}), it follows that

\begin{equation}\label{3dm}
A_\varepsilon \int_M |\nabla_g \eta_\varepsilon|^p
u_\varepsilon^p\; dv_g \leq A_\varepsilon
(A_\varepsilon^{-\frac{1}{p}}
u_\varepsilon(y_\varepsilon)^{\frac{p}{n}})^p
\int_{D_{\varepsilon}} u_\varepsilon^p\; dv_g \leq c
u_\varepsilon(y_\varepsilon)^r (A_\varepsilon^{\frac{1}{p}}
u_\varepsilon(y_\varepsilon)^{-\frac{p}{n}})^n \leq c
m_\varepsilon^{- \frac{p^2}{n}}\; .
\end{equation}

\n Consequently, combining these inequalities with (\ref{3dpt}),
one arrives at

\begin{equation}\label{31d}
A_\varepsilon \int_M |\nabla_g u_\varepsilon|^p
\eta_\varepsilon^p\; dv_g + c d_\varepsilon
A_\varepsilon^{\frac{\tau}{p}} ||u_\varepsilon||_{L^p(M)}^{\tau -
p} \int_M u_\varepsilon^p \eta_\varepsilon^p\; dv_g + c
\frac{\int_M u_\varepsilon^p \eta_\varepsilon^p\; dv_g}{\int_M
u_\varepsilon^p dv_g} \leq c m_\varepsilon^{-\frac{p^2}{n}}\, .
\end{equation}

\n Now, since $p > 1$,the non-sharp Moser
inequality produces

\begin{equation}\label{3ddb}
\int_{B(y_\varepsilon, \frac{1}{2}
A_\varepsilon^{\frac{1}{p}}
u_\varepsilon(y_\varepsilon)^{-\frac{p}{n}})} u_\varepsilon^r\;
dv_g \leq \int_M
(u_\varepsilon \eta_\varepsilon^p)^r\; dv_g \leq c \left(\int_M |\nabla_g u_\varepsilon|^p
\eta_\varepsilon^p\; dv_g\right) \left( \int_M (u_\varepsilon
\eta_\varepsilon^p)^p\; dv_g \right)^{\frac{p}{n}}
\end{equation}

\[
+ c \left(\int_M |\nabla_g \eta_\varepsilon|^p u_\varepsilon^p\;
dv_g\right) \left(\int_M (u_\varepsilon \eta_\varepsilon^p)^p\;
dv_g \right)^{\frac{p}{n}} + c \left(\int_M
(u_\varepsilon \eta_\varepsilon^p)^p\; dv_g\right) \left( \int_M
(u_\varepsilon \eta_\varepsilon^p)^p\; dv_g \right)^{\frac{p}{n}}\, .
\]

\n Thanks to (\ref{3dm}), (\ref{31d}) and since $p > 1$, we then
can estimate each term of the right-hand side of (\ref{3ddb}).
Indeed, we have

\[
\int_M |\nabla_g u_\varepsilon|^p \eta_\varepsilon^p \; dv_g
\left( \int_M (u_\varepsilon \eta_\varepsilon^p)^p\; dv_g
\right)^{\frac{p}{n}} \leq A_\varepsilon \int_M
|\nabla_g u_\varepsilon|^p \eta_\varepsilon^p\; dv_g \left(
\frac{\int_M u_\varepsilon^p \eta_\varepsilon^p\; dv_g}{\int_M
u_\varepsilon^p dv_g} \right)^{\frac{p}{n}}
\leq c m_\varepsilon^{-\frac{p^2}{n}(1 + \frac{p}{n})},
\]

\[
\int_M |\nabla_g \eta_\varepsilon|^p u_\varepsilon^p\; dv_g
\left(\int_M (u_\varepsilon \eta_\varepsilon^p)^p\; dv_g
\right)^{\frac{p}{n}} \leq A_\varepsilon \int_M
|\nabla_g \eta_\varepsilon|^p u_\varepsilon^p\; dv_g \left(
\frac{\int_M u_\varepsilon^p \eta_\varepsilon^p\; dv_g}{\int_M
u_\varepsilon^p dv_g} \right)^{\frac{p}{n}}
\leq c m_\varepsilon^{-\frac{p^2}{n}(1 + \frac{p}{n})} \; .
\]

\n Using the fact

\[
A_\varepsilon \int_M u_\varepsilon^p dv_g \rightarrow 0
\]

\n when $\varepsilon \rightarrow 0$ and because $\tau\leq p$,
we have

\[
A_\varepsilon^{\frac{\tau - p}{p}}
||u_\varepsilon||_{L^p(M)}^{\tau - p} = \left(A_\varepsilon \int_M
u_\varepsilon^p dv_g\right)^{\frac{\tau - p}{p}}
> c
> 0 \; ,
\]

\n for all $\varepsilon$, so

\[
\int_M (u_\varepsilon \eta_\varepsilon^p)^p\; dv_g \left( \int_M
(u_\varepsilon \eta_\varepsilon^p)^p\; dv_g \right)^{\frac{p}{n}} \leq A_\varepsilon \int_M u_\varepsilon^p
\eta_\varepsilon^p\; dv_g \left(\frac{ \int_M u_\varepsilon^p
\eta_\varepsilon^p \; dv_g}{\int_M u_\varepsilon^p dv_g}
\right)^{\frac{p}{n}} \leq c
m_\varepsilon^{-\frac{p^2}{n}(1 + \frac{p}{n})}\ .
\]

\n Replacing these three estimates in (\ref{3ddb}), one gets

\[
\int_{B(y_\varepsilon, \frac{1}{2}
A_\varepsilon^{\frac{1}{p}}
u_\varepsilon(y_\varepsilon)^{-\frac{p}{n}})} u_\varepsilon^r\;
dv_g \leq c m_\varepsilon^{-\frac{p^2}{n}(1
+ \frac{p}{n})}\, ,
\]

\n so that

\[
m_\varepsilon^{\frac{p^2}{n}} \int_{B(y_\varepsilon, \frac{1}{2}
A_\varepsilon^{\frac{1}{p}}
u_\varepsilon(y_\varepsilon)^{-\frac{p}{n}})} u_\varepsilon^r\;
dv_g \leq c m_\varepsilon^{-\frac{p^3}{n^2}}\, .
\]

\n Since $m_\varepsilon \rightarrow \infty$, we obtain that

\[
m_\varepsilon^{\frac{p^2}{n}} \int_{B(y_\varepsilon, \frac{1}{2}
A_\varepsilon^{\frac{1}{p}}
u_\varepsilon(y_\varepsilon)^{\frac{p}{n}})} u_\varepsilon^r\;
dv_g \rightarrow 0\, ,
\]

\n when $\varepsilon \rightarrow 0$.  But this contradicts
(\ref{3dcm}).

\subsection{The final argument in the proof of Theorem \ref{tgno1}} \label{final-contradiction}

In the sequel, we will perform several estimates by using the
$L^r$-concentration and the pointwise estimation. By simplicity,
we consider that the radius of injectivity of $M$ is greater than
one.

Let $\eta \in C^1_0(\R)$ be a cutoff function as in the previous
section and define $\eta_{\varepsilon}(x) =
\eta(d_g(x,x_\varepsilon))$. From the inequality
$M_E(A(p,n))$, we have

\[
\int_{B(0,1)} u_\varepsilon^r \eta_{\varepsilon}^r\; dx
\leq A(p,n) \left( \int_{B(0,1)} |\nabla (u_\varepsilon \eta_{\varepsilon})|^p\; dx\right)
\left(\int_{B(0,1)} u_\varepsilon^p \eta_{\varepsilon}^p\; dx
\right)^{\frac{p}{n}}.
\]

\n Expanding the metric $g$ in normal coordinates around
$x_\varepsilon$ (Cartan expansion), one locally gets

\begin{equation} \label{car}
(1 - c d_g(x,x_\varepsilon)^2)\; dv_g \leq dx \leq (1 + c
d_g(x,x_\varepsilon)^2)\; dv_g
\end{equation}

\n and

\begin{equation} \label{car1}
|\nabla(u_\varepsilon \eta_\varepsilon)|^p \leq
|\nabla_g(u_\varepsilon \eta_\varepsilon)|^p(1 + c
d_g(x,x_\varepsilon)^2) \, .
\end{equation}

\n Thus,

\[
\int_{B(0,1)} u_\varepsilon^r \eta_{\varepsilon}^r\; dx
\]

\[
\leq  \left( A_\varepsilon A(p,n) \int_M
|\nabla_g (u_\varepsilon \eta_{\varepsilon})|^p\; dv_g + c
A_\varepsilon \int_M |\nabla_g
(u_\varepsilon\eta_{\varepsilon})|^p d_g(x,x_\varepsilon)^2\; dv_g
\right) \left(\frac{\int_{B(0,1)} u_\varepsilon^p
\eta_{\varepsilon}^p\; dx}{\int_M u_\varepsilon^p\;
dv_g}\right)^{\frac{p}{n}}\, .
\]

\n Applying then the inequalities

\[
|\nabla_g (u_\varepsilon \eta_{\varepsilon})|^p \leq |\nabla_g
u_\varepsilon|^p \eta_{\varepsilon}^p + c |\eta_{\varepsilon}
\nabla_g u_\varepsilon|^{p - 1} |u_\varepsilon \nabla_g
\eta_{\varepsilon}| + c |u_\varepsilon \nabla_g
\eta_{\varepsilon}|^p \; ,
\]

\n (\ref{3ep}) and (\ref{rm}), we have

\[
A(p,n) \; A_\varepsilon \int_M |\nabla_g
u_\varepsilon|^p\; dv_g  \leq 1 - d_\varepsilon \left( A_\varepsilon
\int_M u_\varepsilon^p\; dv_g \right)^{\frac{\tau}{p}}.
\]

\n One easily checks that

\[
\int_{B(0,1)} u_\varepsilon^r \eta_{\varepsilon}^r\; dx
\]
\begin{equation}\label{3du}
\leq  \left(1 -  d_\varepsilon \left(A_\varepsilon \int_M
u_\varepsilon^p\; dv_g\right)^{\frac{\tau}{p}} + c F_\varepsilon +
c G_\varepsilon + c A_\varepsilon \int_{B(x_\varepsilon,1)
\setminus B(x_\varepsilon,\frac{1}{2})} u_\varepsilon^p \;
dv_g\right)\left(\frac{\int_{B(0,1)} u_\varepsilon^p
\eta_{\varepsilon}^p\; dx}{\int_M u_\varepsilon^p\;
dv_g}\right)^{\frac{p}{n}} \; ,
\end{equation}

\n where

\[
F_\varepsilon = A_\varepsilon \int_M |\nabla_g u_\varepsilon|^p
\eta_{\varepsilon}^p d_g(x,x_\varepsilon)^2\; dv_g \; ,
\]

\n and

\[
G_\varepsilon = A_\varepsilon \int_M |\nabla_g
u_\varepsilon|^{p-1} \eta_{\varepsilon}^{p-1} u_\varepsilon
|\nabla_g \eta_{\varepsilon}|\; dv_g\, .
\]

\n We now estimate $F_\varepsilon$ and $G_\varepsilon$.

\n Applying H\"{o}lder and Young's inequalities, we
obtain

\begin{equation}\label{ge}
G_\varepsilon \leq \kappa \; A_\varepsilon \int_M |\nabla_g
u_\varepsilon|^p  \eta_\varepsilon^p d_g(x,x_\varepsilon)^2 dv_g +
c(\kappa)\; A_\varepsilon \int_{B(x_\varepsilon,1) \setminus
B(x_\varepsilon,\frac{1}{2})} u_\varepsilon^p dv_g  =  \kappa \;
F_\varepsilon + c(\kappa) A_\varepsilon \int_{B(x_\varepsilon,1)}
u_\varepsilon^p dv_g \; ,
\end{equation}

\n where $1 > \kappa > 0$.

Now we will divide into two cases:

{\bf (i) \; $1 < p \leq 2$\;:}

Using H\"{o}lder and Young's inequalities, we obtain

\[
A_\varepsilon \int_M |\nabla_g u_\varepsilon|^{p - 1}
\eta_{\varepsilon}^p u_\varepsilon d_g(x,x_\varepsilon)\; dv_g
\leq \left(A_\varepsilon \int_M |\nabla_g u_\varepsilon|^p
\eta_\varepsilon^p d_g(x,x_\varepsilon)^2 dv_g \right)^{\frac{p -
1}{p}} \left(A_\varepsilon \int_{B(x_\varepsilon,1)}
u_\varepsilon^p dv_g\right)^{\frac{1}{p}}
\]
\begin{equation}\label{3dfq2}
\leq \kappa F_\varepsilon + c(\kappa) A_\varepsilon
\int_{B(x_\varepsilon,1)} u_\varepsilon^p dv_g\; .
\end{equation}

{\bf (ii) $2 < p$\;:}

Using H\"{o}lder and Young's inequalities and as $A_\varepsilon \int_M
|\nabla_g u_\varepsilon|^p dv_g$ is bounded, we obtain

\[
A_\varepsilon \int_M |\nabla_g u_\varepsilon|^{p - 1}
\eta_\varepsilon^p u_\varepsilon d_g(x,x_\varepsilon) \; dv_g \leq
A_\varepsilon \int_M \eta_\varepsilon^{\frac{p}{2}} |\nabla_g
u_\varepsilon|^{\frac{p}{2}} d_g(x,x_\varepsilon) u_\varepsilon
|\nabla_g u_\varepsilon|^{\frac{p - 2}{p}} dv_g
\]
\[
\leq A_\varepsilon \left(\int_M \eta_\varepsilon^p |\nabla_g
u_\varepsilon|^p d_g(x,x_\varepsilon)^2 dv_g\right)^{\frac{1}{2}}
\left(\int_M u_\varepsilon^p dv_g \right)^{\frac{1}{p}} \left(
\int_M |\nabla_g u_\varepsilon|^p dv_g \right)^{\frac{p-2}{2p}}
\]
\begin{equation}\label{efe}
\leq \kappa F_\varepsilon + c(\kappa) \left(A_\varepsilon \int_M u_\varepsilon^p dv_g\right)^{\frac{2}{p}}.
\end{equation}

\n Now taking $u_\varepsilon d_g^2 \eta_\varepsilon^p$ as a test
function in (\ref{3ep}), one easily checks that

\[
F_\varepsilon = A_\varepsilon \int_M |\nabla_g u_\varepsilon|^p
\eta_{\varepsilon}^p d_g(x,x_\varepsilon)^2\; dv_g \leq c
\int_{B(x_\varepsilon, 1)} u_\varepsilon^r
d_g(x,x_\varepsilon)^2\; dv_g + c A_\varepsilon \int_M |\nabla_g
u_\varepsilon|^{p - 1} \eta_{\varepsilon}^p u_\varepsilon
d_g(x,x_\varepsilon)\; dv_g + c G_\varepsilon\, .
\]

\n Therefore, by (\ref{ge}) and (\ref{3dfq2}), we have that

\begin{equation}\label{pequeno}
F_\varepsilon \leq c \int_{B(x_\varepsilon, 1)} u_\varepsilon^r
d_g(x,x_\varepsilon)^2\; dv_g + c(\kappa) A_\varepsilon
\int_{B(x_\varepsilon,1)} u_\varepsilon^p
dv_g \; ,
\end{equation}

\n for $1 < p \leq 2$.

\n On the other hand,  if $2 < p$, we can use the same idea above
together with (\ref{efe}) and we obtain

\begin{equation}\label{grande}
F_\varepsilon \leq c \int_{B(x_\varepsilon, 1)} u_\varepsilon^r
d_g(x,x_\varepsilon)^2\; dv_g + c(\kappa) \left(A_\varepsilon
\int_{B(x_\varepsilon,1)} u_\varepsilon^p
dv_g\right)^{\frac{2}{p}}\; .
\end{equation}

\n Now by (\ref{car}) and the mean value theorem, we obtain

\[
\int_M u_\varepsilon^r \eta_\varepsilon^r\; dv_g
\geq \int_M u_\varepsilon^r
\eta_\varepsilon^r\; dv_g - c \int_M u_\varepsilon^r
\eta_\varepsilon^r d_g(x,x_\varepsilon)^2\; dv_g
\]

\[
\geq 1 - c\int_{M \setminus B(x_\varepsilon,1)} u_\varepsilon^r\;
dv_g - c \int_M u_\varepsilon^r \eta_\varepsilon^r
d_g(x,x_\varepsilon)^2\; dv_g \; ,
\]

\n and by Cartan expansion

\[
\left( \frac{\int_{B(x_\varepsilon,1)} u_\varepsilon^p
\eta_\varepsilon^p\; dx}{\int_M u_\varepsilon^p\; dv_g}
\right)^{\frac{p}{n}} \leq \left( \frac{\int_M
u_\varepsilon^p \eta_\varepsilon^p\; dv_g + c \int_M
u_\varepsilon^p \eta_\varepsilon^p d_g(x,x_\varepsilon)^2\;
dv_g}{\int_M u_\varepsilon^p\; dv_g} \right)^{\frac{p}{n}}
\]

\[
\leq \left( \frac{\int_M u_\varepsilon^p \eta_\varepsilon^p\;
dv_g}{\int_M u_\varepsilon^p\; dv_g} \right)^{\frac{p}{n}} + c \frac{\int_M u_\varepsilon^p
\eta_\varepsilon^p d_g(x,x_\varepsilon)^2\; dv_g}{\int_M
u_\varepsilon^p\; dv_g} \leq 1 + c \frac{\int_M u_\varepsilon^p
\eta_\varepsilon^p d_g(x,x_\varepsilon)^2\; dv_g}{\int_M
u_\varepsilon^p\; dv_g}\, .
\]

\n Using $u_\varepsilon \eta_\varepsilon^p d_g(x_\varepsilon,x)^2$ as the test function in (\ref{3ep}), we also have

\[
\frac{\int_M u_\varepsilon^p \eta_\varepsilon^p
d_g(x,x_\varepsilon)^2\; dv_g}{\int_M u_\varepsilon^p\; dv_g} \leq c \int_{B(x_\varepsilon,1)} u_\varepsilon^r d_g(x_\varepsilon,x)^2 dv_g + F_\varepsilon \; .
\]

\n Replacing these two estimates in (\ref{3du}), one gets

\[
d_\varepsilon \left(A_\varepsilon \int_M u_\varepsilon^p\;
dv_g\right)^{\frac{\tau}{p}}
\]
\begin{equation}\label{final}
\leq c A_\varepsilon \int_M
u_\varepsilon^p dv_g + c F_\varepsilon + c \int_M u_\varepsilon^r
\eta_\varepsilon^r d_g(x,x_\varepsilon)^2 \; dv_g + c \int_{M
\setminus B(x_\varepsilon,1)} u_\varepsilon^r\; dv_g \, .
\end{equation}

\n  {\bf If (iii) $1 < p \leq 2$\;:} using the pointwise
estimate and the definition of $r$, we have

\[
\int_M u_\varepsilon^r \eta_\varepsilon^p d_g(x,x_\varepsilon)^2 dv_g \leq \int_M u_\varepsilon^{r - p} d_g(x,x_\epsilon)^p u_\varepsilon^p dv_g \leq c A_\varepsilon \int_M u_\epsilon^p dv_g \;,
\]

\n and proceeding analogously

\[
\int_{M\setminus B(x_\varepsilon,1)} u_\varepsilon^r dv_g \leq \int_M u_\varepsilon^{r - p} d_g(x,x_\varepsilon)^p u_\varepsilon^p dv_g \leq c A_\varepsilon \int_M u_\varepsilon^p dv_g.
\]

\n {\bf If (iv) $p > 2$\;:} using the pointwise estimate, we have

\[
\int_M u_\varepsilon^r d_g(x,x_\varepsilon)^2 dv_g = \int_M u_\varepsilon^{\frac{2(r - p)}{p}} d_g(x,x_\varepsilon)^2 u_\varepsilon^{\frac{r(p - 2) + 2p}{p}} dv_g \leq c A_\varepsilon^{\frac{2}{p}} \int_M u_\varepsilon^{\frac{r(p - 2) + 2p}{p}} dv_g \; ,
\]

\n and because $p < \frac{r(p - 2) + 2p}{p} < r$, we have by
interpolation that

\[
\int_M u_\varepsilon^r d_g(x,x_\varepsilon)^2 dv_g \leq c \left(A_\varepsilon \int_M u_\varepsilon^p dv_g \right)^{\frac{2}{p}}
\]

\n and proceeding analogously

\[
\int_{M \setminus B(x_\varepsilon,1)} u_\varepsilon^r dv_g \leq \int_M u_\varepsilon^r d_g(x,x_\varepsilon)^2 dv_g \leq c \left(A_\varepsilon \int_M u_\varepsilon^p dv_g \right)^{\frac{2}{p}} \; .
\]

\n Consider $1 < p < 2$. Using (i) and (iii) in (\ref{final}), we
obtain

\begin{equation}\label{f1}
c < |M|^{- \frac{\tau}{n}} - \gamma_\varepsilon \leq
\frac{B_\varepsilon - \gamma_\varepsilon}{c_\varepsilon} =
d_\varepsilon \leq c \left(A_\varepsilon \int_M u_\varepsilon^p
dv_g\right)^{\frac{p - \tau}{p}} \; ,
\end{equation}

\n for some constant $c > 0$.

\n When $2 < p$, we can use (ii) and (iv) in (\ref{final}), and we
have

\begin{equation}\label{f2}
c < |M|^{- \frac{\tau}{n}} - \gamma_\varepsilon \leq \frac{B_\varepsilon - \gamma_\varepsilon}{c_\varepsilon} = d_\varepsilon \leq c \left(A_\alpha \int_M u_\varepsilon^p dv_g \right)^{\frac{2 - \tau}{p}} \; ,
\end{equation}

\n for some constant $c > 0$.

\n If $\tau < \min\{p,2\}$, we have a contradiction when
$\varepsilon \rightarrow 0$. Therefore (A) occurs and ${\cal B} <
\infty$. When $\tau = \min\{p,2\}$, we have by (\ref{f1}) or by
(\ref{f2}) that $d_\varepsilon$ is bounded, \textit{i.e.}, $B_\varepsilon$ is
bounded. This finishes the proof of Theorem 1. \bl

\section{Existence of extremals for the optimal Moser inequality}

Now, let us prove Theorem \ref{extremal}.

In the proof of the Theorem \ref{tgno1}, if ${\cal B} =
M^{-\frac{\tau}{n}}$ we have that the constant functions are
extremal. If ${\cal B} > M^{- \frac{\tau}{n}}$, then we have the
condition (C.2) (see section \ref{EL}) of Theorem 1. As by
hypothesis $\tau < \min\{p,2\}$, we have that the condition
(A) occurs (see section \ref{EL} again), then

\[
\lim_{\varepsilon \rightarrow 0} A_\varepsilon > 0 \; .
\]

\n Using this fact in
(\ref{3ep}), we see that there exists $c > 0$ such that

\[
\int_M|\nabla_g u_\varepsilon|^p dv_g + \left(\int_M
u_\varepsilon^p dv_g\right)^{\frac{\tau}{p}} \leq c \;
\]

\n for all $\varepsilon$. Then, up to a subsequence,
$u_\varepsilon \rightharpoonup u_0$ in $H^{1,p}(M)$. In
particular, since $||u_\varepsilon||_{L^r(M)} = 1$ for all
$\varepsilon$, we have that $||u_0||_{L^r(M)} = 1$, \textit{i.e.}, $u_0$ is
not a null function.

\n From  (\ref{3ep}), we have

\[
\int_M |\nabla_g u_\varepsilon|^{p - 2} \nabla_g u_\varepsilon
\nabla_g h dv_g \leq c \int_M u_\varepsilon^{r - 1} h dv_g \; ,
\]

\n for an arbitrary test function $h \geq 0$. Then, by Moser's iterative scheme for $p < n$ or Morrey's inequality for $p \geq n$ (see \cite{serrin})), we have

\[
\sup_{x \in M} u_\varepsilon \leq c \; ,
\]

\n for all $\varepsilon > 0$ and $p > 1$. Since $u_\varepsilon
\in L^\infty(M)$ for all $\varepsilon$, we can use the Tolksdorf's
regularity in (\ref{3ep}) so that it follows $u_\varepsilon
\rightarrow u_0$ in $C^1(M)$.

\n By (\ref{eq1}) and (\ref{3dha}), the function $v_\varepsilon$
(as obtained in Theorem 1) satisfies

\[
\left(\int_M v_\varepsilon^r dv_g \right)^{\frac{\tau}{p}} \geq
\left( A(p,n)^{\frac{\tau}{p}} \left(\int_M |\nabla_g
v_\varepsilon|^p dv_g \right)^{\frac{\tau}{p}} + (B_\varepsilon -
\gamma_\varepsilon) \left(\int_M v_\varepsilon^p dv_g
\right)^{\frac{\tau}{p}} \right) \left(\int_M v_\varepsilon^p dv_g
\right)^{\frac{\tau}{n}} \; ,
\]

\n and since $u_\varepsilon =
\frac{v_\varepsilon}{||v_\varepsilon||_{L^r(M)}}$, we get

\[
1 \geq \left( A(p,n)^{\frac{\tau}{p}} \left(\int_M |\nabla_g
u_\varepsilon|^p dv_g \right)^{\frac{\tau}{p}} + (B_\varepsilon -
\gamma_\varepsilon) \left(\int_M u_\varepsilon^p dv_g
\right)^{\frac{\tau}{p}} \right) \left(\int_M u_\varepsilon^p dv_g
\right)^{\frac{\tau}{n}} \; .
\]

\n Taking the limit in this inequality, we find

\[
1 \geq \left( A(p,n)^{\frac{\tau}{p}} \left(\int_M |\nabla_g
u_0|^p dv_g \right)^{\frac{\tau}{p}} + {\cal B} \left(\int_M u_0^p
dv_g \right)^{\frac{\tau}{p}} \right) \left(\int_M u_0^q dv_g
\right)^{\frac{\tau}{n}} \; .
\]

\n Then $u_0$ is an extremal function for
$M_R(A(p,n)^{\frac{\tau}{p}},{\cal B})$ and ${\cal B} = B_{opt}$.
\bl

\end{document}